\documentstyle[12pt]{article}

\begin{document}
\newcommand{\ol }{\overline}
\newcommand{\ul }{\underline }
\newcommand{\ra }{\rightarrow }
\newcommand{\lra }{\longrightarrow }
\newcommand{\ga }{\gamma }
\newcommand{\st }{\stackrel }
\title{{\bf The Higher Schur-Multiplicator of Certain Classes of Groups}
\footnote{This research was in part supported by a grant from IPM, Iran.}}
\author{Mohammad Reza R. Moghaddam, B. Mashayekhy and S. Kayvanfar \\
Faculty of Mathematical Sciences, Ferdowsi University of Mashhad,\\
P.O.Box 1159-91775, Mashhad, Iran\\
E-mail: \{fmoghadam, mashaf, skayvanf\}@math.um.ac.ir}
\date{ }
\maketitle
\begin{abstract}
 The paper is devoted to calculating the higher Schur-multiplicator of certain classes of groups with respect to the variety of
 nilpotent groups. Our results somehow generalize the works of M.R.R. Moghaddam (1979), and N.D. Gupta and M.R.R. Moghaddam (1993).\\
\ \ \\
{\it AMS Mathematical Subject Classification}: 20E10, 20F12, 20F18.\\
{\bf Keywords}: Schur-multiplicator, variety of groups, nilpotent product.\\
\end{abstract}

\section{Introduction}

Let ${\cal N}_c$ be the variety of nilpotent groups of class at most $c$ ($c\geq 1$) and $G$ be an arbitrary group with a free presentation
$$ 1\lra R\lra F\stackrel{\pi}{\lra} G\lra 1,$$
where $F$ is a free group of countable rank and $R=\ker{\pi}$. Then the Baer-invariant of $G$ with respect to the variety of nilpotent groups,
${\cal N}_c$, is defined to be
$$ {\cal N}_cM(G)=R\cap \ga_{c+1}/[R,\ _cF ],   $$
where $\ga_{c+1}(F)$ is the ($c+1$)st-term of the lower central series of $F$ and $[R,\ _cF]$ denotes the commutator subgroup $[R,\underbrace {F,F,\ldots ,F}_{c-times}]$.
One may check that ${\cal N}_cM(G)$ is abelian and independent of the choice of the free presentation of $G$ (see [1, 4, 5]).
Determining these Baer-invariants of a given group is known to b every useful for classification of groups into isologism classes (see [2, 4]). In 1992,
Gupta and Moghaddam [3] calculated the Baer-invariants of nilpotent dihedral groups of class $c$ ($c>2$) with respect to the variety of nilpotent groups.
In this paper we calculate the Baer-invariants of the $n^{th}$-nilpotent product of two cyclic groups for $n =2,3,4$, under some conditions.
Our results generalize the works of Moghaddam [5] and Gupta and Moghaddam [3].

\section{Preliminaries and Some Technical Lemmas}

Let $G=G_2\star G_1$ be the free product of two arbitrary groups $G_1$ and $G_2$. Then the $n^{th}$-nilpotent product of $G_1$
and $G_2$ is defined as follows (see also [7]):
$$G_2\stackrel{n}{\star}G_1=\frac{G}{[G_2,G_1,\ _{(n-1)}G]}\ \ for\ n>1.$$
Let $1\lra R_i\lra F_i\stackrel{\pi_i}{\lra} G_i\lra 1$ be a free presentation for $G_i$ $(i=1, 2)$,
where $F_1$ and $F_2$ are free groups, and $R_i=\ker\pi_i$.
We then obtain the following free presentation for the $n^{th}$-nilpotent product $G_2\star G_1$ as follows:
$$1\lra R\lra F\lra G_2\star G_1\lra 1,$$
where $F=F_2\star F_1$ and
$$R=<R_2,R_1,[F_2,F_1,\ _{(n-1)}F]>^F=R_2R_1\prod_{i=1}^{2}[R_i,F][F_2,F_1,\ _{(n-1)}F].$$

Now, in the above set up, assume that
$$Z_r=<x|x^r=1>\ ,\ Z_s=<y|y^s=1>$$
be two cyclic groups of orders $r$ and $s$, respectively, and consider the free presentations for
$Z_r$, $Z_s$ and $Z_r\st{n}{\star} Z_s$, when $F_1=<x>$, $F_2=<y>$, $R_1=<x^r>^{F_1}$, and $R_2=<y^s>^{F_2}$.
Clearly,
$$Z_r\st{n}{\star} Z_s=<x,y|x^r,y^s,\ga_{n+1}(F)>, $$
which is denoted by $G_{(r,s,n)}$. Now put $S=<R_1,R_2>^F$, and $\rho_{c+1}(S)=[S,\ _cF]$ for $c\geq 0$.
Then
$$S=\rho_1(S)\supset \rho_2(S) \supset \rho_3(S) \ldots  $$
is a central series in the free group $F$ . One notes that
$$ S=<R_1,R_2>^F=R_1R_2\prod_{i=1}^{2}[R_i,F]$$
and hence in this case $R=S\ga_{n+1}(F)$.

In this paper we determine the Baer-invariants of the group $G_{(r,s,n)}$
with respect to the variety of nilpotent groups, in different cases, when $c\geq n$.
We keep all the notations throughout the rest of the paper. The following technical lemma is vital in our investigation.

\hspace{-0.65cm}\textbf{Lemma 2.1.}
Let $F$ be a free group and freely generated by $\{x,y\}$. Then the following congruence holds, for all positive integers,
$c\geq 3$, $r\geq 4$, and $a_i\in\{x,y\}$, modulo $\ga_{c+5}(F)$;  
$$ [x^r,y,a_1,\ldots, a_{c-1}]\equiv [x,y,a_1,\ldots, a_{c-1}]^r[x,y,x,a_1,\ldots, a_{c-1}]^{{r \choose 2}}$$
$$ [x,y,a_1,[x,y],a_2,\ldots, a_{c-1}]^{{r \choose 2}}[x,y,a_1,a_2,[x,y,a_1],a_3,\ldots, a_{c-1}]^{{r \choose 2}}$$
$$ [x,y,x,a_1,[x,y],a_2,\ldots, a_{c-1}]^{{r+1 \choose 3}}[x,y,x,a_1,[x,y,x],a_2,\ldots, a_{c-1}]^{{r \choose 3}+{r+1 \choose 3}}$$
$$[x,y,x,x,a_1,\ldots, a_{c-1}]^{{r \choose 3}} [x,y,x,[x,y],a_1,\ldots, a_{c-1}]^{{r \choose 3}}$$
$$[x,y,x,x,x,a_1,\ldots, a_{c-1}]^{{r \choose 4}},$$
where ${r \choose k}$ denotes $r!/k!(r-k)!$.\\
\textbf{Proof.} We use double inductions, and so we only need to show the validity of
the congruent for $c=3$, since the induction on $c$ completes the proof. Hence assume $c=3$, the proof goes by induction on $r$ ($r\geq 4$).
By expanding the left hand side for $r=4$ and working modulo  
$\ga_8(F)$, we obtain the formula rather easily. Now assume the result holds for $r-1$ ($r\geq 5$), we obtain the following congruence,
modulo $\ga_8(F)$
$$ [x^{r-1},y,a_1,a_2]\equiv[[x^{r-1},y][x^{r-1},y,x][x,y],a_1,a_2]$$
$$\equiv [x^{r-1},y,a_1,a_2][x^{r-1},y,a_1,a_2,[x,y,a_1]]$$
$$ [x^{r-1},y,a_1,[x,y],a_2][x^{r-1},y,a_1,[x,y,a_1]^{r-1},a_2]$$
$$[x^{r-1},y,x,a_1,a_2][x^{r-1},y,x,a_1,[x,y],a_2][x,y,a_1,a_2]$$
Clearly, modulo $\ga_8(F)$
$$ [x^{r-1},y,a_1,a_2,[x,y,a_1]]\equiv[x,y,a_1,a_2,[x,y,a_1]]^{r-1}$$
$$ [x^{r-1},y,a_1,[x,y],a_2]\equiv[x,y,a_1,[x,y],a_2]^{r-1}$$
$$ [x^{r-1},y,a_1,[x,y,x]^{r-1},a_2]\equiv[x,y,a_1,[x,y,x],a_2]^{(r-1)^ 2}$$
$$ [x^{r-1},y,x,a_1,[x,y],a_2]\equiv[x,y,x,a_1,[x,y],a_2]^{r-1}$$
Now, using the induction formula for $r-1$ to the following commutators
$$ [x^{r-1},y,a_1,a_2]\ and\ [x^{r-1},y,x,a_1,a_2]$$
we obtain the desired formula.$\Box$

Now we are able to prove our key result, which shortens the proof of the main theorems in the next section.

\hspace{-0.65cm}\textbf{Proposition 2.2.}
Let $Z_r$ and $Z_s$ be two cyclic groups of orders $r$ and $s$, respectively. Let $F_1/R_1\cong Z_r$,
$F_2/R_2\cong Z_s$, $F=F_1\star F_2$, $S=<R_1,R_2>^F$ and $\rho_{c+1}(S)=[S,\ _{c+1}F]$ $(c\geq 1)$. Then\\
$(i)$ $\ga_{c+2}(F)\rho_{c+1}(S)/\ga_{c+2}(F)=dZ\oplus \ldots \oplus dZ$ $(r(c+1)\ copies)$, for all 
nonnegative integers $r$ and $s$, and $d=(r,s)$.\\
$(ii)$ $\ga_{c+j}(F)\rho_{c+1}(S)/\ga_{c+j}(F)=dZ\oplus \ldots \oplus dZ$ $(\sum_{i=1}^{i=j-1}r(c+i)\ copies)$,
where $j=3,4,5$ and for $j=3$ , then $r$ and $s$ must be odd, and for $j=4$ or $5$ ,
then r and s should not be divisible by $2$ and $3$. Also $r(c+i)$ is the number
of basic commutators of weight $c+i$ on two letters.\\
\textbf{Proof.} $(i)$ Modulo $\ga_{c+2}(F)$, and using induction argument the following congruences hold:
$$ [x^r,a_1,\ldots ,a_c]\equiv[x,a_1,\ldots ,a_c]^r\ ;$$
$$ [y^s,a_1,\ldots ,a_c]\equiv[y,a_1,\ldots ,a_c]^s\ ,$$
for all $a_1,\ldots ,a_c\in F$.
Clearly these commutators are elements of $\rho_{c+1}(S)$. We note that $\ga_{c+2}(F)\rho_{c+1}(S)/\ga_{c+2}(F)$  
is generated by all the $r^{th}$ and $s^{th}$-powers,
and hence $(r,s)=d^{th}$-powers of the basic commutators of weight $c+1$. Thus
the property of being abelian gives the result.\\
$(ii)$ If $j=3$, then we have the following exact sequence
$$ 1\lra \frac{\ga_{c+3}(F)\rho_{c+2}(S)}{\ga_{c+3}(F)}\lra \frac{\ga_{c+3}(F)\rho_{c+1}(S)}{\ga_{c+3}(F)}\lra \frac{\ga_{c+3}(F)\rho_{c+1}(S)}{\ga_{c+3}(F)\rho_{c+2}(S)}\lra 1.$$
By part $(i)$, $\ga_{c+3}(F)\rho_{c+2}(S)/\ga_{c+3}(F)$ is an abelian group generated by all the
$d^{th}$-powers of the basic commutators of weight $c+2$ on two letters. Now we have
the following congruence:
$$ [x^r,a_1,\ldots ,a_c]\equiv[x,a_1,\ldots ,a_c]^r[x,a_1,x,a_2,\ldots ,a_c]^{{r \choose 2}}\ (mod\ \ga_{c+3}(F)). $$
Since $r$ is odd, we obtain
$$ [x,a_1,x,a_2,\ldots ,a_c]^{{r \choose 2}}\equiv [x^r, a_1,x,a_2,\ldots ,a_c]^{\frac{r-1}{2}}\ (mod\ \ga_{c+3}(F)). $$
which is in $\rho_{c+2}(S)$. Hence modulo $\ga_{c+3}(F)\rho_{c+2}(S)$,
$$ [x^r,a_1,\ldots ,a_c]\equiv[x,a_1,\ldots ,a_c]^r\ ,$$
$$ [y^s,a_1,\ldots ,a_c]\equiv[y,a_1,\ldots ,a_c]^s\ .$$
Therefore $\ga_{c+3}(F)\rho_{c+1}(S)/\ga_{c+3}(F)\rho_{c+2}(S)$ is freely generated by all the $d^{th}$-powers of the basic commutators of weight $c+1$. 
Hence the above split extension gives the result for $j=3$.

If $j=4$, then we have the following exact sequences:
$$ 1\ra E=\frac{\ga_{c+4}(F)\rho_{c+2}(S)}{\ga_{c+4}(F)\rho_{c+3}(S)}\ra B=\frac{\ga_{c+4}(F)\rho_{c+1}(S)}{\ga_{c+4}(F)\rho_{c+3}(S)}\ra D=\frac{\ga_{c+4}(F)\rho_{c+1}(S)}{\ga_{c+4}(F)\rho_{c+2}(S)}\ra 1, $$
$$ 1\ra C=\frac{\ga_{c+4}(F)\rho_{c+3}(S)}{\ga_{c+4}(F)}\ra A=\frac{\ga_{c+4}(F)\rho_{c+1}(S)}{\ga_{c+4}(F)}\ra B=\frac{\ga_{c+4}(F)\rho_{c+1}(S)}{\ga_{c+4}(F)\rho_{c+3}(S)}\ra 1, $$
which are also split extensions. Using the pro of of part $(i)$ and the case $j=3$, we obtain
$$ E\cong dZ\oplus \ldots \oplus dZ\ \ \ \ (r(c+2)-copies),$$
$$ C\cong dZ\oplus \ldots \oplus dZ\ \ \ \ (r(c+3)-copies).$$
Now for every generator $[x^r,y,a_1,\ldots ,a_{c-1}]$ or $[y^s,x,a_1,\ldots ,a_{c-1}]$ of $\rho_{c+1}(S)$
and using Lemma 2.1 we have, modulo $\ga_{c+4}(F)$,
$$ [x^r,y,a_1,\ldots, a_{c-1}]\equiv [x,y,a_1,\ldots, a_{c-1}]^r[x,y,x,a_1,\ldots, a_{c-1}]^{{r \choose 2}}$$
$$ [x,y,a_1,[x,y],a_2,\ldots, a_{c-1}]^{{r \choose 2}}[x,y,x,x,a_1,\ldots, a_{c-1}]^{{r \choose 3}}. $$
Since $2$ and $3$ do not divide $r$ , it implies that
$$ [x^r,y,a_1,\ldots, a_{c-1}]\equiv [x,y,a_1,\ldots, a_{c-1}]^r\ \ (mod \ga_{c+4}(F)\rho_{c+2}(S)).$$
Similarly
$$ [y^s,x,a_1,\ldots, a_{c-1}]\equiv [y,x,a_1,\ldots, a_{c-1}]^s\ \ (mod \ga_{c+4}(F)\rho_{c+2}(S)).$$
Thus
$$ D\cong dZ\oplus \ldots \oplus dZ\ \ \ \ (r(c+1)-copies).$$
Now the above exact sequences give the result.

Finally, for $j=5$ we consider the following split exact sequences
$$ 1\ra E_i=\frac{\ga_{c+5}(F)\rho_{c+i}(S)}{\ga_{c+5}(F)\rho_{c+i+1}(S)}\ra B_i=\frac{\ga_{c+5}(F)\rho_{c+1}(S)}{\ga_{c+5}(F)\rho_{c+i+1}(S)}$$
$$\ra D_i=\frac{\ga_{c+5}(F)\rho_{c+1}(S)}{\ga_{c+5}(F)\rho_{c+i}(S)}\ra 1,\ \ for\ i=2\ and\ 3. $$
$$ 1\ra C=\frac{\ga_{c+5}(F)\rho_{c+4}(S)}{\ga_{c+5}(F)}\ra A=\frac{\ga_{c+5}(F)\rho_{c+1}(S)}{\ga_{c+5}(F)}\ra B_2=\frac{\ga_{c+5}(F)\rho_{c+1}(S)}{\ga_{c+5}(F)\rho_{c+4}(S)}\ra 1. $$
Following the previous procedure we conclude that
$$ D_2\cong dZ\oplus \ldots \oplus dZ\ \ \ \ (r(c+1)-copies),$$
$$ E_i\cong dZ\oplus \ldots \oplus dZ\ \ \ \ (r(c+2)-copies)(i=2,3),$$
$$ C\cong dZ\oplus \ldots \oplus dZ\ \ \ \ (r(c+4)-copies).$$
Then the above exact sequences guarantee the result.$\Box$

\section{The Main Results}
Using the results in the previous section we are able to calculate the higher Schur-
multiplicators ${\cal N}_cM(G)$, with respect to the variety of nilpotent groups of
class at most $c$ $(c\geq n)$. By the above notations,
$$ {\cal N}_cM(G_{(r,s,n)})=\frac{R\cap\ga_{c+1}(F)}{[R,\ _cF]}=\frac{S\ga_{n+1}(F)\cap\ga_{c+1}(F)}{[S\ga_{n+1}(F),\ _cF]}  $$
$$=\frac{\ga_{c+1}(F)}{\ga_{c+n+1}(F)\rho_{c+1}(S)}\cong \frac{\ga_{c+1}(F)/\ga_{c+n+1}(F)}{\ga_{c+n+1}(F)\rho_{c+1}(S)/\ga_{c+n+1}(F)}.$$
Clearly for $c\geq n$, $\ga_{c+1}(F)/\ga_{c+n+1}(F)$ is free abelian group of rank $\sum_{i=1}^{n}r(c+i)$,
where $r(c+i)$ denotes the rank of the lower central factor group $\ga_{c+i}(F)/\ga_{c+i+1}(F)$.  
One notes that the main problem is to determine the structure of the factor group  
$\ga_{c+n+1}(F)\rho_{c+1}(S)/\ga_{c+n+1}(F)$. Now, if $n=1$, then by Proposition
2.2 (i) we obtain the following theorem.

\hspace{-0.65cm}\textbf{Theorem 3.1.}
Let $r$ and $s$ be arbitrary positive integers. Then for any $c\geq 1$ and $(r,s)=d$,
$$ {\cal N}_cM(G_{(r,s,1)})\cong Z_d\oplus \ldots \oplus Z_d\ \ \ \ (r(c+1)-copies).$$

\hspace{-0.65cm}\textbf{Remark 3.2.}
The above theorem gives the complete structure of ${\cal N}_cM(Z_r\times Z_s)$
for all $c\geq 1$. This generalizes Theorem 3.3 of Moghaddam [5].

The first two authors [6] in 1997 presented an explicit formula for the higher
Schur-multiplicator of a finite abelian group with respect to the variety of nilpotent groups of class at most 
$c\geq 1$, ${\cal N}_c$, as follows:

\hspace{-0.65cm}\textbf{Theorem 3.3.} (M.R.R. Moghaddam and B. Mashayekhy (1997)) 
Let $G=Z_{n_1}\oplus Z_{n_2}\oplus \ldots \oplus Z_{n_k}$ be a finite abelian group, where 
$n_{i+1}\mid n_i$ for all $1\leq i\leq k-1$.
Then, for all $c\geq 1$, the higher Schur-multiplicator of $G$ is
$${\cal N}_cM(G)\cong Z_{n_2}^{(b_2)}\oplus Z_{n_3}^{(b_3-b_2)}\oplus \ldots \oplus Z_{n_k}^{(b_k-b_{k-1})},$$
where $b_i$ is the number of basic commutators of weight $c+1$ on $i$ letters, and $Z_{n}^{(m)}$
denotes the direct sum of $m$ copies of the cyclic group $Z_n$.\\
\textbf{Proof.} See Theorem 2.4 of [6].

\hspace{-0.65cm}\textbf{Theorem 3.4.}
If $(r,s)=1$, then for any $n\geq 1$ and $c\geq 1$,
$$ {\cal N}_cM(G_{(r,s,n)})=1.$$
\textbf{Proof.} Using the commutator identities and induction argument we obtain the following congruences, modulo $\ga_{n+1}(F)$, 
where $F$ is the free group on $\{x,y\}$.   
$$ [a_1,\ldots ,a_{i-1},x^r,a_{i+1},\ldots ,a_n]\equiv[a_1,\ldots ,a_{i-1},x,a_{i+1},\ldots ,a_n]^r\ ,$$
$$ [a_1,\ldots ,a_{i-1},y^s,a_{i+1},\ldots ,a_n]\equiv[a_1,\ldots ,a_{i-1},y,a_{i+1},\ldots ,a_n]^s\ ,$$
for all $a_1,\ldots ,a_n\in \{x,y\}$.
Hence all the $r^{th}$ and $s^{th}$-powers of the basic commutators of weight $n$ are trivial in the following
$$G_{(r,s,n)}\cong ,x,y|x^r,y^s,\ga_{n+1}(F)>.$$
Hence $G_{(r,s,n)}\cong ,x,y|x^r,y^s,\ga_{n}(F)>.$ Using the above procedure for the
commutators of weight less than n , and after finite number of steps we deduce that
$$G_{(r,s,n)}\cong ,x,y|x^r,y^s,\ga_{2}(F)>\cong Z_r\times Z_s.$$
Now either of Theorems 3.3 or 3.1 gives the result.$\Box$

Using the notations at the beginning of this section and Proposition 2.2, we can prove the following theorem.

\hspace{-0.65cm}\textbf{Theorem 3.5.}

$(i)$ For any odd integers $r$ and $s$, and all $c\geq 2$
$$ {\cal N}_cM(G_{(r,s,2)})\cong Z_d\oplus \ldots \oplus Z_d\ \ \ \ (\sum_{i=1}^{2}r(c+i)-copies).$$

$(ii)$ For all non-negative integers $r$ and $s$, which are not divisible by $2$ and $3$,
then
$$ {\cal N}_cM(G_{(r,s,3)})\cong Z_d\oplus \ldots \oplus Z_d\ \ \ \ (\sum_{i=1}^{3}r(c+i)-copies),$$
where $c\geq 3$ and
$$ {\cal N}_cM(G_{(r,s,4)})\cong Z_d\oplus \ldots \oplus Z_d\ \ \ \ (\sum_{i=1}^{4}r(c+i)-copies),$$
where $c\geq 4$.

\hspace{-0.65cm}\textbf{Remark 3.6.}
The above theorem determines the complete structure of $n^{th}$-nilpotent product of two cyclic groups, when $n\leq 4$. 
These results generalize Gupta and Moghaddam [3].

Finally, we remark that if one can construct a similar identity congruence
as in Lemma 2.1, one may obtain the same results for $n\geq 5$. This certainly
involves a very complicated commutators manipulations.

\end{document}